\documentclass[a4paper,reqno]{amsart}
\usepackage{amsmath,amsfonts}
\usepackage{amssymb}
\usepackage{amsthm}
\usepackage{graphicx}
\usepackage{tabularx}
\usepackage{array}
\usepackage{physics}
\usepackage{tikz}
\usepackage[colorlinks,linkcolor=cyan,citecolor=cyan]{hyperref}
\usepackage{cleveref}

\crefname{thm}{Theorem}{Theorems}
\crefname{lem}{Lemma}{Lemmas}
\crefname{prop}{Proposition}{Propositions}
\crefname{cor}{Corollary}{Corollaries}
\usepackage{mathtools}
\usepackage{physics}
\usepackage{algorithm}
\usepackage{algorithmic}

\usepackage[square,comma,sort&compress,numbers]{natbib} %% Sort of References
\usepackage{mathrsfs} %% Swash letter
\usepackage{booktabs} %% tables with three lines

\newtheorem{thm}{Theorem}%[section]
\newtheorem{cor}[thm]{Corollary}
\newtheorem{lem}[thm]{Lemma}

\theoremstyle{definition}

\newtheorem{rem}[thm]{Remark}

\newcommand{\SET}[2]{\left\{#1\;\middle|\; #2\right\}}

\DeclareMathOperator{\SPAN}{span}
\DeclareMathOperator{\GL}{GL}

\begin{document}
\title{Banach-Mazur Distance from $\ell_p^3$ to $\ell_\infty^3$}

\date{}

\author{Longzhen Zhang}

\email{zhanglongzhen0404@163.com}

\author{Lingxu Meng}

\email{menglingxu@nuc.edu.cn}

\author{Senlin Wu}

\email{wusenlin@nuc.edu.cn}

\address{College of Mathematics, North University of China, 030051 Taiyuan China}

\begin{abstract}
  The maximum of the Banach-Mazur distance $d_{BM}^M(X,\ell_\infty^n)$, where
  $X$ ranges over the set of all $n$-dimensional real Banach spaces, is
  difficult to compute. In fact, it is already not easy to get the maximum of
  $d_{BM}^M(\ell_p^n,\ell_\infty^n)$ for all $p\in [1,\infty]$. We prove that
  $d_{BM}^M(\ell_p^3,\ell_\infty^3)\leq 9/5,~\forall p\in[1,\infty]$. As an
  application, the following result related to Borsuk's partition problem in
  Banach spaces is obtained: any subset $A$ of $\ell_p^3$ having diameter $1$ is
  the union of $8$ subsets of $A$ whose diameters are at most $0.9$.
\end{abstract}



\maketitle

\section{Introduction}
The (multiplicative) \emph{Banach-Mazur distance} between two isomorphic Banach
spaces $X$ and $Y$ is defined as
\begin{displaymath}
  d_{BM}^{M}(X,Y)=\inf\SET{\norm{T}\cdot\norm{T^{-1}}}{T\text{~ is an isomorphism from
      $X$ onto $Y$}}.
\end{displaymath}
It is well known that
\begin{displaymath}
  d_{BM}^M(X,Y)\leq d_{BM}^M(X,Z)\cdot d_{BM}^M(Z,Y),
\end{displaymath}
where $X$, $Y$, and $Z$ are isomorphic Banach spaces, see, e.g., \cite{Jaegermann1989}.

A compact convex subset of $\mathbb{R}^n$ having interior points is called a
\emph{convex body}. Let $\mathcal{K}^{n}$ be the set of all convex bodies in
$\mathbb{R}^{n}$ and $\mathcal{C}^{n}$ be the set of convex bodies that are
symmetric with respect to the \emph{origin} $o$ of $\mathbb{R}^n$. Let
$\mathcal{A}^{n}$ be the set of all nonsingular affine transformations on
$\mathbb{R}^{n}$. The Banach-Mazur distance between $K,L\in \mathcal{K}^{n}$ is
defined by
\begin{displaymath}  
  d_{BM}^{M}(K,L)=\inf\SET{\gamma \ge 1 }{ \exists T\in \mathcal{A}^n,~x\in
  \mathbb{R}^n,~\text{s.t.} ~ T(L)\subseteq K\subseteq \gamma T(L)+x}.
\end{displaymath}
The infimum can be attained. When $K,L\in \mathcal{C}^n$, one can verify that
\begin{displaymath}
  d_{BM}^{M}(K,L)=\inf\SET{\gamma \ge 1 }{\exists T\in \mathcal{T}^n,~\text{s.t.}~T(L)\subseteq K\subseteq \gamma T(L)},
\end{displaymath}
where $\mathcal{T}^n$ is the set of all nonsingular linear transformations on
$\mathbb{R}^{n}$. Denote by $B_X$ the unit ball of an $n$-dimensional Banach
space $X=(\mathbb{R}^n,\norm{\cdot})$. We have
$d_{BM}^{M}(X,Y)=d_{BM}^{M}(B_{X},B_{Y})$, which connects the Banach-Mazur
distance between finite dimensional Banach spaces with the Banach-Mazur distance
between two convex bodies (cf. e.g., \cite[p. 15, p.
47]{Fabian-Habala-Hajek-Montesinos-Zizler2011}) and provides a link between
Banach space theory and convex geometry. It is generally difficult to calculate the
exact value of the Banach-Mazur distance between convex bodies (or isomorphic
Banach spaces).

Denote by $\ell_{p}^{n}$ the space $(\mathbb{R}^n,\norm{\cdot}_p)$,
where the $p$-norm $\norm{\cdot}_p$ is given by
\begin{displaymath}
  \norm{(\alpha_{1},\cdots,\alpha_{n})}_{p}=\qty(\sum\limits_{i\in
    [n]}{|\alpha_{i}|^{p}})^{\frac{1}{p}},~\forall p\in[1,\infty),
\end{displaymath}
and
\begin{displaymath}
  \norm{(\alpha_{1},\cdots,\alpha_{n})}_{\infty}=\max\limits_{i\in [n]}{|\alpha_{i}|}.
\end{displaymath}
Here we used the shorthand notation $[n]:=\SET{i\in \mathbb{Z}^{+}}{1\le i\le
  n}$. Denote by $B_p^n$ the unit ball of $\ell_{p}^{n}$. Clearly,
$B_{\infty}^{n}=[-1,1]^{n}$. We have the following classical result:
\begin{thm}[{cf. \cite[Proposition 37.6]{Jaegermann1989}}]
  \label{thm:estimations}
  Let $n$ be a positive integer and $1\le p,q \le \infty$. 
  \begin{enumerate}
  \item[(i)] If $1\le p\le q\le 2$ or $2\le p\le q\le \infty$, then $d_{BM}^{M}(\ell_{p}^{n},\ell_{q}^{n})=n^{1/p-1/q}$.
  \item[(ii)] If $1\le p<2<q \le \infty$, then $\gamma n^{\alpha} \le
    d_{BM}^{M}(\ell^{n}_{p},\ell_{q}^{n})\le \eta n^{\alpha}$, where $\alpha=\max\qty{\flatfrac{1}{p}-\flatfrac{1}{2},\flatfrac{1}{2}-\flatfrac{1}{q}}$,
    and $\gamma$, $\eta$ are universal constants. If $n=2^k$
    \emph{(}$k\in\mathbb{N}$\emph{)}, then $\eta=1$.
  \end{enumerate}  
\end{thm}
From \Cref{thm:estimations}, it follows that
$d_{BM}^{M}(\ell^{n}_{p},\ell^{n}_{\infty})=n^{1/p},~\forall p\in [2,\infty]$.
In general, it is difficult to get the exact value of
$d_{BM}^{M}(\ell^{n}_{p},\ell^{n}_{\infty})$ for $p\in [1,2)$. The case when $n=2$ is an
exception. Since $\ell_1^2$ and $\ell_\infty^2$ are isometric, 
\begin{displaymath}
  d_{BM}^{M}(\ell^{2}_{p},\ell^{2}_{\infty})=d_{BM}^{M}(\ell^{2}_{p},\ell^{2}_{1})=2^{1-1/p},~\forall
  p\in  [1,2).
\end{displaymath}
When $n=2^k$ for some $k\in\mathbb{N}$, we have
$d_{BM}^{M}(\ell_{p}^{n},\ell_{\infty}^{n})\leq \sqrt{n},~\forall
p\in[1,\infty]$. In particular, we have
\begin{equation}
  \label{eq:4-d}
  d_{BM}^{M}(\ell_{p}^{4},\ell_{\infty}^{4})\leq 2,~\forall p\in[1,\infty].
\end{equation}
F. Xue \cite{Xue-2018-MR3868093} provided explicit upper bounds of
$d_{BM}^M(\ell_1^n,\ell_\infty^n)$ for $n\in \qty{3,4,5,6,7,8}$. and showed
that
\begin{displaymath}
  \alpha\sqrt{n}\leq d_{BM}^M(\ell_1^n,\ell_\infty^n)\leq
  (\sqrt2+1)\sqrt{n},~\forall n\in \mathbb{Z}^+,
\end{displaymath}
where $\alpha$ is an absolute constant (cf. \cite[Theorem 1.5]{Xue-2018-MR3868093}).

When $n=3$, Y. Lian and S. Wu \cite{Lian-Wu-2021-MR4261748} proved that
\begin{displaymath}
  d_{BM}^{M}(\ell_{p}^{3},\ell_{\infty}^{3}) \le \dfrac{\sqrt{18\cdot
      19}}{10},~\forall p\in [1,2].
\end{displaymath}
In this paper, we improve this result as follows:
\begin{thm}
  \label{thm:main-result}
  We have
  \begin{equation}
    \label{eq:3-d}
    d_{BM}^{M}(\ell_{p}^{3},\ell_{\infty}^{3}) \le \dfrac{9}{5},~\forall  p\in[1,\infty].
  \end{equation}
\end{thm}
Most likely, the estimations \eqref{eq:4-d} and \eqref{eq:3-d}  are both tight.
By \Cref{thm:estimations}, \Cref{thm:main-result}, and \cite[Theorem
2]{Lian-Wu-2021-MR4261748}, we have the following improvment of \cite[Theorem
16]{Lian-Wu-2021-MR4261748}:
\begin{cor}
  \label{cor:Borsuk}
  For each $p\in [1,\infty]$, any set $A$ of $\ell_p^3$ having diameter
  $1$ is the union of $8$ subsets of $A$ whose diameters are at most $0.9$.
\end{cor}
This result is closely related to Borsuk's partition problem in finite
dimensional Banach spaces, see \cite{Zong-2021-MR4338243,Lian-Wu-2021-MR4261748}
for more details. For Borsuk's problem in $\ell_2^3$, Tolmachev et al.
\cite{Tolmachev-Protasov-Voronov-2022-MR4441253} proved that, if the diameter of
$A\subseteq \ell_2^3$ is $1$, then $A$ can be partitioned into four subsets
whose diameters are at most 0.966. Note that, a closed ball in $\ell_\infty^3$
cannot be split into 7 subsets having smaller diameters. Therefore we cannot
replace $8$ with a positive integer $m\leq 7$ and obtain a result similar to
\Cref{cor:Borsuk}.

% In \Cref{sec:optimization-problem}, we show that the Banach-Mazur distance from
% an $n$-dimensional Banach space $(\mathbb{R}^n,\norm{\cdot})$ to $\ell_\infty^n$
% is the optimal value of an optimization problem. Upper bounds of
% $d_{BM}^M(\ell_p^3,\ell_\infty^3)$ for some $p\in [1,2]$ are obtained via an
% inexact algorithm. \Cref{thm:main-result} are proved in \Cref{sec:proof}.

\section{Banach-Mazur distance to $\ell_\infty^n$}
\label{sec:optimization-problem}

Denote by $\GL_n(\mathbb{R})$ the set of all nonsingular $n\times n$ matrices of
real numbers.
For $K,L\in \mathcal{C}^n$ and $A\in \GL_n(\mathbb{R})$, set
\begin{gather*}
  \gamma_1(K,L;A)=\inf\SET{\gamma}{\gamma>0\text{ and }A(L)\subseteq \gamma K},\\
  \gamma_2(K,L;A)=\sup\SET{\gamma}{\gamma>0\text{ and }\gamma K\subseteq A(L)}.
\end{gather*}
Here we identify a member of $\GL_n(\mathbb{R})$ with the corresponding
nonsingular linear transformation. Since both $K$ and $A(L)$ contain the origin
$o$ in their interior, $\gamma_1(K,L;A)$ and $\gamma_2(K,L;A)$ are well-defined
and are positive. Moreover, since $K$ and $A(L)$ are both compact, $\inf$ and
$\sup$ in the definitions above can be replaced with $\min$ and $\max$,
respectively.

Inspired by the proof of \cite[Lemma 14]{Lian-Wu-2021-MR4261748}, we have
\Cref{lem:distance-symmetric-bodies} and \Cref{lem:calculate-gamma}.

\begin{lem}
  \label{lem:distance-symmetric-bodies}
  For $K,L\in\mathcal{C}^n$, 
  \begin{displaymath}
    d_{BM}^M(K,L)=\min\SET{\frac{\gamma_1(K,L;A)}{\gamma_2(K,L;A)}}{A\in \GL_n(\mathbb{R})}.
  \end{displaymath}
\end{lem}
\begin{proof}
  Let $B$ be an arbitrary element of $\GL_n(\mathbb{R})$. By the definitions of
  $\gamma_1(K,L;B)$ and $\gamma_2(K,L;B)$, we have $\gamma_2(K,L;B)K\subseteq B(L)\subseteq \gamma_1(K,L;B)K$, or equivalently,
  \begin{displaymath}
    \qty(\frac{1}{\gamma_1(K,L;B)}B)(L)\subseteq K\subseteq \frac{\gamma_1(K,L;B)}{\gamma_2(K,L;B)}\qty(\frac{1}{\gamma_1(K,L;B)}B)(L).
  \end{displaymath}
  It follows that
  \begin{displaymath}
    d_{BM}^M(K,L)\leq \frac{\gamma_1(K,L;B)}{\gamma_2(K,L;B)}.
  \end{displaymath}
  Hence
  \begin{displaymath}
    d_{BM}^M(K,L)\leq \inf\SET{\frac{\gamma_1(K,L;A)}{\gamma_2(K,L;A)}}{A\in \GL_n(\mathbb{R})}.
  \end{displaymath}
  Conversely, there exists $A_0\in \GL_n(\mathbb{R})$ such that $A_0(L)\subseteq
  K\subseteq d_{BM}^M(K,L)A_0(L)$. Then $\gamma_1(K,L;A_0)\leq 1$ and
  $\gamma_2(K,L;A_0)\geq \qty(d_{BM}^M(K,L))^{-1}$. Hence
  \begin{displaymath}
    d_{BM}^M(K,L)\geq \frac{\gamma_1(K,L;A_0)}{\gamma_2(K,L;A_0)}\geq \inf\SET{\frac{\gamma_1(K,L;A)}{\gamma_2(K,L;A)}}{A\in \GL_n(\mathbb{R})}.
  \end{displaymath}
  This completes the proof.
\end{proof}

Let $(\mathbb{R}^n,\norm{\cdot})$ be a Banach space and let
$(\mathbb{R}^n,\norm{\cdot}_*)$ be its dual. Each $y\in \mathbb{R}^n$
defines a linear functional $f$ on $(\mathbb{R}^n,\norm{\cdot})$ by
$f(x)=y^T\cdot x$. When $y\neq o$, we have
\begin{displaymath}
  \norm{f}_*=\qty(d(o,\SET{x\in \mathbb{R}^n}{f(x)=1}))^{-1},
\end{displaymath}
where $d(\cdot,\cdot)$ is the distance on $\mathbb{R}^n$ induced by $\norm{\cdot}$.
\begin{lem}
  \label{lem:calculate-gamma}
  Let $X=(\mathbb{R}^n,\norm{\cdot})$ be a Banach space and $A=(a_{ij})_{n\times
    n}\in \GL_n(\mathbb{R})$. Denote by $A_{ij}$ the cofactor of $a_{ij}$, and
  set
  \begin{displaymath}
    x_i=(a_{1i},a_{2i},\dots,a_{ni})^\intercal,\quad y_j=(A_{1j},A_{2j},\dots,A_{nj})^\intercal,~\forall i,j\in[n]. 
  \end{displaymath}
Then
  \begin{equation}
    \label{eq:gamma-1}
    \gamma_1(B_X,B_\infty^n;A)=\max\SET{\norm{\sum\limits_{i\in[n]}\sigma_ix_i}}{\sigma_i\in\qty{-1,1},~\forall
      i\in[n]},
  \end{equation}
  \begin{equation}
    \label{eq:gamma-2}
    \gamma_2(B_X,B_\infty^n;A)=\min\SET{\frac{|\det A|}{\norm{y_i}_*}}{i\in[n]}.
  \end{equation}
  In particular,
  \begin{align*}
    &d_{BM}^M(X,\ell_\infty^n)\\
    =&\min\limits_{A\in \GL_n(\mathbb{R})}\max\SET{(|\det
      A|)^{-1}\norm{y_i}_*\norm{\sum\limits_{j\in[n]}\sigma_jx_j}}{i\in[n],~\sigma_j\in\qty{-1,1},~\forall
    j\in [n]}.
  \end{align*}
\end{lem}
\begin{proof}
  Evidently, $A(B_\infty^n)=A([-1,1]^n)$ is a convex polytope with
  \begin{displaymath}
    \SET{\sum\limits_{i\in[n]}\sigma_ix_i}{\sigma_i\in\qty{-1,1},~\forall i\in[n]}
  \end{displaymath}
  as the set of vertices. Moreover, $A(B_\infty^n)$ is contained in $\gamma B_X$
  if and only if every vertex of $A(B_\infty^n)$ is contained in $\gamma B_X$.
  Thus \eqref{eq:gamma-1} holds.

  For each $i\in[n]$, let $H_i$ be the hyperplane passing through $x_i$ and
  parallel to the hyperplane spanned by $\SET{x_j}{j\in [n]\setminus\qty{i}}$.
  We easily verify that, $\pm H_1,\dots,\pm H_n$ are the bounding hyperplanes of
  $A(B_\infty^n)$. For each $i\in[n]$, the null space of the linear functional
  $f_i$ defined by $f_i(x)=y_i^\intercal\cdot x$ is precisely
  $\SPAN\SET{x_j}{j\in[n]\setminus\qty{i}}$. Thus
  \begin{displaymath}
    H_i=\SET{x\in \mathbb{R}^n}{y_i^\intercal\cdot x=y_i^\intercal\cdot x_i}=\SET{x\in\mathbb{R}^n}{\frac{y_i^\intercal}{\det A}\cdot x=1}.
  \end{displaymath}
  It follows that
  \begin{displaymath}
    d(o,H_i)=\qty(\norm{\frac{y_i}{\det A}}_*)^{-1}=\frac{|\det A|}{\norm{y_i}_*}.
  \end{displaymath}
  Assume that $\gamma>0$. Then $\gamma B_X\subseteq A(B_\infty^n)$ if and only
  if $\gamma\leq \min \SET{d(o,H_i)}{i\in[n]}$. Hence the equality \eqref{eq:gamma-2}
  follows.
\end{proof}
\begin{rem}
  Clearly, $\gamma_2(B_X,B_\infty^n;A)$ is the reciprocal of the operator norm
  $\norm{A^{-1}}$ of $A^{-1}$. We can also deduce \eqref{eq:gamma-2} using the
  fact that $\norm{A^{-1}}$ equals the operator norm of its adjoint (cf. e.g.,
  \cite[Lemma 9.1]{Clason-2020-MR4182425}).
  % \begin{displaymath}
  %   \norm{A^{-1}}=\norm{\qty(A^{-1})^*}.
  % \end{displaymath}
\end{rem}

\begin{rem}
  Set
  \begin{displaymath}
    R_\infty^n=\max\SET{d_{BM}^M(X,\ell_\infty^n)}{X\text{ is an
        $n$-dimensional Banach space}}.
  \end{displaymath}
  It is shown in \cite{Giannopoulos-1995-MR1353450} that there exists a
  universal constant $c>0$ such that $R_\infty^n\leq cn^{5/6}$. S. Taschuk
  \cite{Taschuk-2011-MR2794363} proved that, for $n\geq 3$,
  \begin{equation}
    \label{eq:distance-low-dimension}
    R_\infty^n \le  \sqrt{n^2-2n+2+\cfrac{2}{\sqrt{n+2}-1}}.
  \end{equation}
  P. Youssef \cite{Youssef-2014-MR3164527} showed that $R_\infty^n\leq
  (2n)^{5/6}$, which is better than the estimation in
  \eqref{eq:distance-low-dimension} when $n\geq 22$. \Cref{lem:calculate-gamma}
  provides a way for estimating $R_\infty^n$ when $n$ is small.
\end{rem}

% If $1\le p,~q\le \infty$ and $\flatfrac{1}{p}+\flatfrac{1}{q}=1$, then we say
% that $p$ and $q$ are \emph{conjugate}.

\section{Banach-Mazur distance from $\ell_p^3$ to $\ell_\infty^3$}
\label{sec:proof}

Assume that $A=(a_{ij})_{3\times 3}\in \GL_3(\mathbb{R})$ and $A_{ij}$ is the
cofactor of $a_{ij},~\forall i,j\in[3]$. Let $x_1,x_2,x_3$ be the
column vectors of $A$ and set $y_i=(A_{1i},A_{2i},A_{3i})^\intercal,~\forall
i\in[3]$. For $p\in [1,\infty]$, put
\begin{equation}
  \label{eq:formula}
  g_p(A)=\frac{1}{|\det A|}\max\SET{\norm{y_i}_q\norm{x_1+\sigma_2x_2+\sigma_3x_3}_p}{i\in[3],~\sigma_1,\sigma_2\in\qty{-1,1}},
\end{equation}
where $q$ is the conjugate of $p$. Set
$d(p)=d_{BM}^M(\ell_p^3,\ell_\infty^3),~\forall p\in [1,\infty]$. By
\Cref{lem:calculate-gamma}, $d(p)$ is the optimal value of the optimization
problem
\begin{equation}
  \label{eq:optimization}
  \min\limits_{A\in \GL_3(\mathbb{R})}\quad g_p(A).
\end{equation}
By \eqref{eq:distance-low-dimension}, $d(p)\leq \flatfrac{\sqrt{2\left(
      \sqrt{5}+11 \right)}}{2}\approx 2.572553$. Put
\begin{displaymath}
  \mathcal{J}=\SET{A\in \GL_3(\mathbb{R})}{g_p(A)\leq \frac{\sqrt{2\left( \sqrt{5}+11 \right)}}{2}}.
\end{displaymath}
Then \eqref{eq:optimization} is equivalent to the optimization problem
\begin{displaymath}
  \min\limits_{A\in \mathcal{J}}\quad g_p(A).
\end{displaymath}

We use the Nelder-Mead simplex algorithm (cf.
\cite{Rios-Sahinidis-2013-MR3070154,Lagarias-Reeds-Wright-Wright-1999-MR1662563})
to find a local minimum of $g_p(A)$ starting from some $A\in \mathcal{J}$, and
apply a particle swarm algorithm (cf. \cite{Gazi-Passino-2011-MR3235758}) to
process a global search. Numerical experiments yield estimations for upper
bounds of $d(p)$, see \Cref{tab:estimations}. When $p=2$, the estimation in
\Cref{tab:estimations} is very close to $\sqrt{3}$, which is the exact value of
$d(2)$.

\begin{table}
  \centering
  \begin{tabular}{c|c|c|c|c|c|c}
    \toprule
    $p$ & 1 & 1.2 & 1.4 & 1.6 & 1.8 & 2 \\
    \midrule
    upper bound of $d(p)$  & 1.8000 & 1.71533 & 1.67744 & 1.67601 & 1.69732 & 1.73205\\
    \bottomrule
  \end{tabular}  
  \caption{Several estimations of $d(p)$}
  \label{tab:estimations}
\end{table}

\begin{lem}
  \label{lem:first-portion}
  For $p\in [1,1.7]$, $d(p)\le \flatfrac{9}{5}$.
\end{lem}
\begin{proof}
  By the proof of \cite[Lemma 14]{Lian-Wu-2021-MR4261748}, we have
  \begin{gather}
    d(1)\leq \frac{\norm{(1,4,1)}_1\norm{(3,1,3)}_\infty}{10}=\frac{9}{5},\nonumber\\
    d(p)\le \frac{1}{10}(4^p+2)^{\frac{1}{p}}\cdot (2\cdot 3^{\frac{p}{p-1}}+1)^{\frac{p-1}{p}}, \forall p\in (1,2].\label{eq:dp}
  \end{gather}
  Thus we only need to consider the case when $p\in (1,1.7]$. Set
  \begin{displaymath}
    f(p)=\ln(2+4^p)+(p-1)\cdot \ln(2\cdot 3^{\frac{p}{p-1}}+1),~\forall p\in(1,2]. 
  \end{displaymath}
  Then $(4^p+2)^{\frac{1}{p}}\cdot (2\cdot
  3^{\frac{p}{p-1}}+1)^{\frac{p-1}{p}}=e^{\frac{f(p)}{p}}$. For $p\in (1,2]$, put
 $r(p)=\flatfrac{f(p)}{p}$ and $w(p)=pf'(p)-f(p)$. We have
  $r'(p)=\flatfrac{w(p)}{p^2}$ and $w'(p)=pf''(p)$,
  where
  \begin{gather*}
    f'(p)=\frac{4^p}{2+4^p}\cdot \ln 4+\ln (2+3^{-\frac{p}{p-1}})+\ln
    3+\frac{\ln 3}{(p-1)(2\cdot 3^{\frac{p}{p-1}}+1)},\\
    f''(p)=\frac{2\cdot \ln^2 4\cdot 4^p}{(2+4^p)^2}+\frac{\ln 3}{(p-1)^3}\cdot \frac{2 \cdot \ln 3 \cdot 3^{\frac{p}{p-1}}}{(2\cdot 3^{\frac{p}{p-1}}+1)^2}.
  \end{gather*}  
  Obviously, $\lim\limits_{p\to 1^+}f(p)=\ln 18$ and $\lim\limits_{p\to
    1^+}f'(p)=\frac{2}{3}\cdot \ln 4+\ln 6$. Therefore,
  \begin{equation}
    \label{eq:w(1)}
    \lim\limits_{p\to 1^+}w(p)=\frac{2}{3}\cdot\ln 4-\ln 3<0.
  \end{equation}
  Moreover,
  \begin{equation}
    \label{eq:w(2)}
    w(2)=\frac{16}{9}\ln 4+2\ln 19-\frac{36}{19}\ln 3-\ln (18\cdot 19)>0. %\tag{4.4}
  \end{equation}
  Since $f''(p)$ is positive on $(1,2]$, $w(p)$ is strictly increasing on
  $(1,2]$. By \eqref{eq:w(1)} and \eqref{eq:w(2)}, there exists a unique $p_0\in
  (1,2)$ satisfying $w(p_0)=0$. Therefore, $r'(p)\le 0$ for $p\in (1,p_0]$ and
  $r'(p)>0$ for $p\in (p_0,2]$. Hence $r(p)$ decreases on $(1,p_0]$ and
  increases on $(p_0,2]$. Since
  \begin{displaymath}
    2.8904\approx\ln 18=\lim\limits_{p\to 1^+}r(p)>r(1.7)\approx 2.8864,
  \end{displaymath}
  we have $r(p)\le \lim\limits_{p\to 1^+}r(p)=\ln 18,~\forall p\in (1,1.7]$. By \eqref{eq:dp},
  \begin{displaymath}
    d(p)\le \dfrac{e^{r(p)}}{10} \le \dfrac{e^{\ln 18}}{10} =
    \dfrac{9}{5},~\forall p\in (1,1.7]. \qedhere
  \end{displaymath}
\end{proof}

Now we are ready to prove \Cref{thm:main-result}.
\begin{proof}[Proof of \Cref{thm:main-result}]
  By \Cref{thm:estimations} and \Cref{lem:first-portion}, we only need to
  consider the case when $p\in [1.7,2]$. Set
  \begin{displaymath}
    A_1=
    \begin{pmatrix}
      13 & -24 & 24\\
      -24 & 13 & 24\\
      24 & 24 & 13
    \end{pmatrix}\qqtext{and}A_2=
    \begin{pmatrix}
      9 & -17 & 17\\
      -17 & 9 & 17\\
      17 & 17 & 9
    \end{pmatrix}.
  \end{displaymath}
  Using \eqref{eq:formula}, we get $d(1.7) \le g_{1.7}(A_1) \le 1.6967$ and
  $d(1.8) \le g_{1.8}(A_2) \le 1.7033$. By \Cref{thm:estimations},
  \begin{gather*}
    d_{BM}^{M}\qty(\ell_{1.7}^{3},\ell_{p}^{3}) =3^{1/1.7-1/p}\le
    3^{1/1.7-1/1.8}\le 1.0366,~\forall p\in[1.7,1.8],\\
    d_{BM}^{M}\left( \ell_{1.8}^{3},\ell_{p}^{3} \right) =3^{1/1.8-1/p}\le
    3^{1/1.8-1/1.9}\le 1.0327,~\forall p\in [1.8,1.9],\\
    d_{BM}^{M}\left( \ell_{2}^{3},\ell_{p}^{3} \right) =3^{1/p-1/2}\le
    3^{1/1.9-1/2}\le 1.0294,~\forall p\in[1.9,2].
  \end{gather*}
  It follows that
  \begin{gather*}
    d_{BM}^M(\ell_p^3,\ell_\infty^3)\leq
    d_{BM}^M(\ell_{1.7}^3,\ell_\infty^3)\cdot
    d_{BM}^{M}\qty(\ell_{1.7}^{3},\ell_{p}^{3})<\frac{9}{5},~\forall p\in
    [1.7,1.8],\\
    d_{BM}^M(\ell_p^3,\ell_\infty^3)\leq
    d_{BM}^M(\ell_{1.8}^3,\ell_\infty^3)\cdot
    d_{BM}^{M}\qty(\ell_{1.8}^{3},\ell_{p}^{3})<\frac{9}{5},~\forall p\in
    [1.8,1.9],\\
    d_{BM}^M(\ell_p^3,\ell_\infty^3)\leq
    d_{BM}^M(\ell_{2}^3,\ell_\infty^3)\cdot
    d_{BM}^{M}\qty(\ell_{2}^{3},\ell_{p}^{3})<\frac{9}{5},~\forall p\in [1.9,2].
  \end{gather*}
  Thus $d(p) \le \flatfrac{9}{5}, ~\forall p\in [1.7,2]$. This completes the proof.
\end{proof}

\section{Acknowledgement}
The authors are supported by the National Natural Science Foundation of China
(grant numbers 12071444 and 12001500), the Natural Science Foundation of Shanxi
Province of China (grant numbers 201901D111141 and 202103021223191), and the
Scientific and Technological Innovation Programs of Higher Education
Institutions in Shanxi (grant number 2020L0290).

%\bibliographystyle{amsplain}
%\bibliography{../xbib/references}
\def\polhk#1{\setbox0=\hbox{#1}{\ooalign{\hidewidth
  \lower1.5ex\hbox{`}\hidewidth\crcr\unhbox0}}} \def\cprime{$'$}
  \def\polhk#1{\setbox0=\hbox{#1}{\ooalign{\hidewidth
  \lower1.5ex\hbox{`}\hidewidth\crcr\unhbox0}}} \def\cprime{$'$}
\providecommand{\bysame}{\leavevmode\hbox to3em{\hrulefill}\thinspace}
\providecommand{\MR}{\relax\ifhmode\unskip\space\fi MR }
\providecommand{\MRhref}[2]{%
  \href{http://www.ams.org/mathscinet-getitem?mr=#1}{#2}
}
\providecommand{\href}[2]{#2}

\end{document}